\documentclass[11pt, letterpaper]{article}
\usepackage[utf8]{inputenc}
\usepackage{amsmath}
\usepackage{amsfonts}
\usepackage{mathtools}
\usepackage{amssymb}
\usepackage{amsthm}
\usepackage{graphicx}
\usepackage[percent]{overpic}  
\usepackage{subcaption}
\usepackage{sidecap}
\usepackage[colorinlistoftodos]{todonotes}
\usepackage{wrapfig}
\usepackage[numbers]{natbib}
\usepackage{multirow}

\DeclareMathOperator{\diag}{diag}

\usepackage[top=1.0in, bottom=1.0in, left=1.0in, right=1.0in]{geometry}

\newtheorem{thm}{Theorem}

\providecommand{\keywords}[1]
{
  \small	
  \textbf{\textit{Keywords---}} #1
}

\title{Universal differential equations for optimal control problems and its application on cancer therapy}
\author{Wenjing Zhang, Wandi Ding, Huaiping Zhu}
\date{\today}

\begin{document}

\maketitle
\begin{abstract}
This paper highlights a parallel between the forward–backward sweeping method for optimal control and deep learning training procedures. We reformulate a classical optimal control problem, constrained by a differential equation system, into an optimization framework that uses neural networks to represent control variables. We demonstrate that this deep learning method adheres to Pontryagin’s Maximum Principle and mitigates numerical instabilities by employing backward propagation instead of a backward sweep for the adjoint equations. As a case study, we solve an optimal control problem to find the optimal combination of immunotherapy and chemotherapy. Our approach holds significant potential across various fields, including epidemiology, ecological modeling, engineering, and financial mathematics, where optimal control under complex dynamic constraints is crucial.
\end{abstract}
 \hspace{10pt}

\section{Deep Learning for Optimal Control: A Pontryagin’s Maximum Principle Perspective}
Classical optimal control problems are typically tackled using direct methods \cite{betts2010practical}, Pontryagin’s Maximum Principle \cite{pontryagin2018mathematical}, or Hamilton-Jacobi-Bellman equations  \cite{bellman1966dynamic}, which provide strong theoretical underpinnings but can struggle with high-dimensional or highly sensitive systems \cite{lenhart2007optimal,betts2010practical,yin2023optimal,caillau2023algorithmic}.
Incorporating modern machine learning techniques into classical optimal control frameworks holds the potential to revolutionize the way we address complex, high-dimensional systems. By leveraging deep learning alongside Pontryagin’s Maximum Principle (PMP), we establish a robust methodology that circumvents common numerical instabilities while maintaining rigorous theoretical guarantees.

Pontryagin’s Maximum Principle (PMP)  is a foundational tool for solving deterministic optimal control problems, where system dynamics are described by ordinary, partial, or difference equations. These problems are governed by specified state equations, control variables, and objective functionals. PMP provides necessary conditions for optimality by defining a Hamiltonian that encompasses state, control, and adjoint variables. Solving the resulting optimality system typically involves a forward-backward sweep \cite{lenhart2007optimal}, which iterates as follows: (1) solve the state equations forward in time using an initial guess for the control; (2) solve the adjoint equations backward in time using boundary or transversality conditions from the optimization of the objective functional; and (3) update the control via PMP’s optimality condition.  

A well-known challenge in this setting is numerical instability when integrating the adjoint equations backward in time, particularly in highly sensitive systems (such as those near bifurcations). Even small computational errors can rapidly escalate, making solutions from backward-in-time integrators (e.g., the backward Euler method) unreliable \cite{chen2018analog,mcasey2012convergence}.


Interestingly, the forward-backward structure in optimal control closely parallels the training procedure in deep learning. By representing the control as a neural network and treating the objective functional as the loss function, the forward pass of a neural network—where inputs propagate through successive layers—mirrors solving the state equations forward in time and accumulating the necessary state constraints. Conversely, the backward pass, which computes gradients of the loss function with respect to each parameter, parallels solving adjoint equations backward in time. The training process iteratively adjusts neural network parameters until changes in the objective functional become negligible.


This idea is closely related to  universal differential equations (UDEs) \cite{chen2018neural,rackauckas2017adaptive,rackauckas2020universal}, which combine known physical laws (represented by differential equations) with unknown mechanisms modeled through universal approximators \cite{hornik1989multilayer,cybenko1989approximation}, such as neural networks. In this context, the neural network can represent the controls, while the state equations capture established physical laws. As a result, constraints and controls seamlessly integrate into a unified framework.


To address the numerical instability in solving the adjoint equation backward in time, backpropagation in deep learning circumvents many of the issues often encountered in backward-in-time integration of adjoint equations. It provides exact gradients that power standard optimizers, including stochastic gradient descent (SGD) \cite{amari1993backpropagation}, momentum-based methods \cite{diakonikolas2021generalized} (such as SGD with momentum \cite{ramezani2024generalization}), adaptive methods \cite{goodfellow2016deep} (such as AdaGrad \cite{duchi2011adaptive}, RMSProp \cite{hinton2012neural}, and Adam \cite{kingma2014adam}), and approximate second-order techniques \cite{lecun2002efficient}, such as conjugate gradients \cite{moller1993scaled,le2011optimization}, the Broyden-Fletcher-Goldfarb-Shanno (BFGS) algorithm \cite{fletcher2000practical}, and limited memory BFGS (L-BFGS) \cite{liu1989limited}.


In addition, common training strategies—such as mini-batch training \cite{lecun2012efficient}, data augmentation \cite{shorten2021text}, and various regularization approaches (for example, dropout \cite{srivastava2014dropout})—as well as learning-rate scheduling \cite{bengio2012practical}, further enhance performance and generalization. By modeling the control as a neural network and optimizing it via backpropagation, practitioners can often circumvent the instability issues inherent in some forward-backward PMP formulations, thereby leveraging the robustness of gradient-based optimization methods.


Scientific machine learning (such as Physics-Informed Neural Networks (PINNs)) blends data-driven techniques with mechanistic models, enabling physically consistent, data-efficient, and interpretable models suitable for a range of forward and inverse problems. While representing a notable advance, PINNs can be computationally costly, numerically tricky, and reliant on good prior physics knowledge \cite{farea2024understanding,huang2022applications}. To overcome these constraints, UDEs integrate universal approximators \cite{chen2018neural} directly into differential equations, offering efficient solvers and greater numerical stability.
This synergy opens promising avenues in optimal control, preserving classical rigor while leveraging the flexibility of modern machine learning approaches. Although previous studies \cite{yin2023optimal,boltz_symbolicoptimalcontrol} have applied UDEs, a gap remains in establishing that an optimal control problem formulated as a UDE comply with Pontryagin’s Maximum Principle (PMP).

In this work, we briefly review the definition of deep neural networks and their training in Sections \ref{deepNN} and \ref{TrainNN}, then define the optimal control problem within a UDE framework in Section \ref{defineUDE}. Next, we show that optimizing the control problem by training neural networks ensures the necessary conditions for PMP in Section \ref{proofOCasDNN}. In Section \ref{Application}, we also provide an application in which we design an optimal treatment strategy for combination immunotherapy for cancer. This model features complex bifurcations and rich dynamical behaviors \cite{zhang2024bifurcations}, which can cause numerical instabilities when using the forward-backward sweep method. By optimizing our UDE formulation, we demonstrate a stable and efficient training approach for solving the optimal control problem. 
Finally, we conclude with a discussion in Section \ref{concl}.

\section{Background}
Deep learning is a subfield of machine learning that utilizes deep neural networks (DNNs), neural networks with multiple layers, to model and approximate complex functions. The theoretical foundation of deep learning lies in the universal approximation theorem \cite{hornik1989multilayer,cybenko1989approximation}, which states that a neural network with sufficient neurons and appropriate training can approximate a wide variety of functions with arbitrary accuracy. Practically, this implies that a properly configured DNN with adequate data can learn nearly any mapping between inputs and outputs.

However, training deep neural networks poses several notable challenges \cite{cheng2018model}. DNNs, particularly those with millions or billions of parameters, require substantial computational resources, making them expensive and difficult to train. Additionally, deep models trained on limited datasets are prone to overfitting, where they memorize specific training examples rather than generalize to new data. Excessive parameters can also introduce redundancy, increasing computational costs without necessarily improving performance. 

UDEs offer a promising approach to addressing these challenges by blending deep learning techniques with scientific modeling. The UDE framework integrates mechanistic models, which are typically differential equations derived from physical principles, with neural networks that capture unknown or difficult-to-model dynamics. This hybrid structure reduces computational inefficiency and reduces the risks of overfitting and redundant parameters inherent to purely data-driven methods \cite{rackauckas2017adaptive,rackauckas2020universal}.

Given these advantages, UDEs provide an efficient framework for solving optimal control problems by deep learning. To set the stage for formulating an optimal control problem within the UDE framework and subsequently establishing that this formulation aligns with Pontryagin’s Maximum Principle (PMP), we now provide a concise overview of deep neural networks, their training methodologies, and the fundamental concepts underlying UDEs.

\subsection{Deep learning and Deep neural network}\label{deepNN}
DNNs consist of an input layer, multiple hidden layers, and an output layer \cite{sze2017efficient}. Each layer is made up of units, or neurons, which are connected to neurons in the previous and subsequent layers. Each connection is assigned a weight, a parameter that determines the influence of one neuron on another. 
We now provide a general formulation of DNN as a foundation for the rest of the texts.

Let $ x \in \mathbb{R}^{n_1} $ denote the input vector, where $ n_1 $ is the dimensionality of the input data. A network with $ L $ layers has layers 1 being the input layer and layer $L$ being the output layer.
Each layer $ l $, for $l=2,3, \dots, L$, consists of $n_l$ neurons that perform weighted sums of the inputs from the previous layer, followed by a non-linear activation function.
To calculate weighted sums, each layer $ l $, for $l=2,3, \dots, L$, has the matrix of weights, also named as weight matrices and denoted as $ W^{(l)} \in \mathbb{R}^{n_{l} \times n_{l-1}} $, and the vector of bias, also named as the bias vector and denoted as $ b^{(l)} \in \mathbb{R}^{n_{l}} $, where $ n_{l} $ represents the number of neurons in layer $ l $, and $ n_{l-1} $ represents the number of neurons in the preceding layer.
Then, we define the pre-activation $ z^{(l)} $ at layer $ l $, which is the linear transformation  $z^{(l)} = W^{(l)} a^{(l-1)} + b^{(l)}$, where $z^{(l)} =\{z_j^{(l)}\}_{j=1}^{j=n_l} \in \mathbb{R}^{n_{l}}$. 
Each element $z_j^{(l)}$ denotes the weighted input for neuron $j$ at layer $l$ with $j=1, \dots, n_l$  \cite{higham2019deep}. 
The action of the network from neuron $j$ at layer $l$ results the output $a^{(l)}_j$. 
During forward propagation, the relation between the input layer through each hidden layer to the output layer then be written as
\begin{equation}
\begin{array}{rll}
a^{(1)} &= x & \text{input layer}, \\[2.0ex]
z^{(l)} &= W^{(l)} a^{(l-1)} + b^{(l)} & \text{pre-activation at layer }l  = 2, 3, \dots, L,  \\[2.0ex]
a^{(l)} &= \{a_j^{(l)}\}_{j=1}^{j=n_l} =  \sigma .\left( z^{(l)}  \right)  & \text{output or activation at layer }l = 2, 3, \dots, L, 
\end{array}
    \label{eqn_NN1}
\end{equation} 
where $\sigma$ is referred to as the nonlinear activation function and $l=2,3, \dots, L$.
The $.$ denote the element-wise operation, that is, if $x,y \in \mathcal{R}^n$, then $x.y \in \mathcal{R}^n$ is defined by $(x.y)_i=x_iy_i$.  
Often, the final layer $a^{L}$ is an output of a final activation function depending on the task.
For regression, the final activation is often taken as the same activation function for the hidden layers.
For the purpose of classification, a softmax function, $a_j^{(l)} = e^{z_j^{(L)}}/\sum_{i} e^{z_i^{(L)}}$, may be applied to produce class probabilities.
Mathematically, the network can be represented as a series of vector-matrix operations followed by a nonlinear activation function at each neuron.
The process that input $x$ is passed through multiple layers to produce an output prediction $z^{(L)}(x)$ is called forward propagation (or forward pass) \cite{zell1994simulation,svozil1997introduction}.

\subsection{Training a neural network}\label{TrainNN}
Training a neural network involves adjusting its parameters (weights and biases) to minimize the difference between the network’s predictions and the actual data. From a data-fitting perspective, this process includes three main stages. First, the forward pass propagates inputs through the neural network from the input layer to the output layer, resulting in predictions. Next, the network’s outputs are evaluated using a predefined loss function to quantify how well they match the actual data. Finally, during the backward pass (backpropagation), the network parameters are iteratively adjusted to minimize the loss.

Mathematically, consider a dataset containing $N$ training points $ \{ x^{(i)}\}_{i=1}^N\subset \mathbb{R}^{n_1}$, each with corresponding desired outputs $\{ y(x^{(i)})\}_{i=1}^N \in \mathbb{R}^{n_L}$.
Forward pass involves a chain of linear transformations (or vector-matrix operations) and activation, then calculates the output from the final layer $\{ a^{(L)}(x^{(i)})\}_{i=1}^N$.

Next, to measure how accurately the neural network maps inputs $\{ x^{(i)}\}_{i=1}^N$ to desired outputs $\{ y(x^{(i)})\}_{i=1}^N$, a cost function is used. A common choice for regression is the mean squared error (MSE):
\begin{equation}
\mathcal{L}(y, a^{(L)}) = \frac{1}{N} \sum_{i=1}^N \mathcal{L}_i(y, a^{(L)}), \quad\text{where} \quad \mathcal{L}_i(y, a^{(L)}) = \| y(x^{(i)}) - a^{(L)}(x^{(i)})\|^2.
    \label{eqn_NN2}
\end{equation} 
We gather all entries from each weight matrix  $W^{(l)}=[w^{(l)}_{jk}] \in \mathbb{R}^{n_l\times n_l}$, and bias vector, $b^{(l)}=[b^{(l)}_{j}] \in \mathbb{R}^{n_l}$, for $l = 2, 3, \dots, L$, into a single parameter vector 
$$
p = \bigl(w_{jk}^{(l)},\, b_{j}^{(l)}\bigr)_{j,k=1,\dots,n_l;\; l=2,\dots,L} \in \mathbb{R}^s,
$$
where $s$ is the total number of parameters in the entire network. 
For a fixed training point $x^{(i)}$ with desired output $y(x^{(i)})$, $i = 1,\dots,N$, the cost function term $\mathcal{L}_i(y, a^{(L)})$ in \eqref{eqn_NN2} becomes $\mathcal{L}_i(p)$, highlighting its dependence on the parameter vector. 

Minimizing the cost function $\mathcal{L}(p): \mathbb{R}^s \rightarrow \mathbb{R}$ involves iterative adjustments of $p$. For a small parameter update $\Delta p$, we approximate using a Taylor series expansion:
\begin{equation}
\mathcal{L}(p+\Delta p) \approx \mathcal{L}(p) + \nabla \mathcal{L}(p)^T\,\Delta p,
    \label{eqn_NN3}
\end{equation}
where
\begin{equation}
 \nabla \mathcal{L}(p) = \left(\frac{\partial \mathcal{L}}{\partial p_r}\right)_{r=1,..,s}
      \label{eqn_NN3b}
\end{equation}
is the gradient of the cost function with respect to the network parameters $p$. 
By the Cauchy-Schwarz inequality, to effectively reduce the cost, updates typically move in the opposite direction of the gradient $- \nabla \mathcal{L}(p)$, leading to the gradient descent update rule:
\begin{equation}
p \quad \rightarrow \quad p \; - \eta \nabla \mathcal{L}(p)^T,
    \label{eqn_NN4}
\end{equation}
where the small stepsize $\eta$ is referred to as the learning rate and $\nabla \mathcal{L}(p) = \frac{1}{N} \sum_{i=1}^N  \nabla \mathcal{L}_i(p)$.
Computing the gradient $\nabla \mathcal{L}(p)$  can be expensive for large datasets, so efficient optimization variants such as stochastic gradient descent (SGD) \cite{amari1993backpropagation}, mini-batch gradient descent \cite{khirirat2017mini}, and adaptive methods like Adam \cite{kingma2014adam,kochenderfer2019algorithms} are widely used.

We randomly choose an initial guess for $p$ and iterate via \eqref{eqn_NN4} until either a convergence threshold is satisfied or we run out of computational budget \cite{higham2019deep}. For instance, a necessary condition \cite{nocedal1999numerical} for a parameter vector $p^{\ast}$ to be a local minimizer of the loss function $\mathcal{L}(p)$ is that

\begin{equation}
 \nabla \mathcal{L}(p^{\ast}) = \left(\frac{\partial \mathcal{L}}{\partial w^{(l)^{\ast}}_{jk}},\,\frac{\partial \mathcal{L}}{\partial b^{(l)^{\ast}}_j}\right)_{j,k=1,\dots,n_l;\; l=2,\dots,L} =0.
    \label{eqn_NN5}
\end{equation}

The backward pass computes the gradient $\nabla \mathcal{L}(p)$. We begin at the output layer $L$ by taking the derivative of the cost  $L$ with respect to the final activations  $ a^{(L)} $, often referred to as the error:
\begin{equation}
    \delta^{(L)} = \frac{\partial \mathcal{L}}{\partial z^{(L)}} = \frac{\partial \mathcal{L}}{\partial a^{(L)}}  .  \sigma'(z^{(L)}),   
    \label{eqn_NN6}
\end{equation}
where $ \sigma'(z^{(L)}) $ is the derivative of the activation function for layer $ L $.
Then, we propagate the error backward through the network.
Specifically, for each previous layer $ l = L-1, \dots, 1 $, we calculate the error term
\begin{equation}
     \delta^{(l)} = \left( W^{(l+1)} \right)^T \delta^{(l+1)}  .  \sigma'(z^{(l)})= D^{(l)}\left( W^{(l+1)} \right)^T D^{(l+1)}\left( W^{(l+2)} \right)^T \cdots D^{(L-1)}\left( W^{(L)} \right)^T D^{(L)}\left( \frac{\partial \mathcal{L}}{\partial a^{(L)}} \right)^T,
    \label{eqn_NN7}
\end{equation}
where, the superscript $^T$ denotes matrix transpose and $D^{(l)}=\diag (\sigma'(z^{(l)}_i) \in \mathbb{R}^{n_l \times n_l}$ \cite{higham2019deep}.
Finally, we calculate the gradients  $\nabla \mathcal{L}(p) = \bigl(\frac{\partial \mathcal{L}}{\partial w^{(l)}_{jk}},\,\frac{\partial \mathcal{L}}{\partial b^{(l)_j}}\bigr)_{j,k=1,\dots,n_l;\; l=2,\dots,L}$, where
\begin{equation}
\frac{\partial \mathcal{L}}{\partial w^{(l)}_{jk}} = \left[\delta^{(l)}  .  \left( a^{(l-1)} \right)^T \right]_{jk}\quad\text{and}\quad
\frac{\partial \mathcal{L}}{\partial b^{(l)_j}} = \left[\delta^{(l)} \right]_j \quad \text{for}\quad l = L-1, \dots, 1.
    \label{eqn_NN8}
\end{equation}

\subsection{Universal differential equations}\label{defineUDE}
In an UDE, the system dynamics are described by
\begin{equation}
\frac{d\hat{y}}{dt}=g(\hat{y},t,q, u(\hat{y},t,p)),
    \label{eqn_NN10}
\end{equation}
where $g$ is a function of the state variables $\hat{y} \in \mathbb{R}^n$, time $t$, parameter set $q$ representing known parameters in the physical model, and $u(\hat{y},t,p)$ is defined by a neural network with parameters $p$. 
One may write $ u(\hat{y},t,p)=a^{(L)}(\hat{y},t,p)$ with output dimension $n_L$.
Embedding a universal approximator $u$ into function $g$ allows us to learn a component of the model from data while preserving the known physical laws that govern the system.

Forward pass involves solving the UDE in \eqref{eqn_NN10} for given parameters $p$ to obtain $\hat{y}(t,p)$.
Suppose there are $N$ time points $\{t^{(i)}\}_{i=1}^N$ with observed data $\{y(t^{(i)})\}_{i=1}^N$.
The solutions of the UDE model \eqref{eqn_NN10} are $\hat{y}(t^{(i)},p)$.
A mean squared error cost is then defined as 
\begin{equation}
\mathcal{L}(p) = \frac{1}{N} \sum_{i=1}^N \| y(t^{(i)}) - \hat{y}(t^{(i)},p)\|^2
    \label{eqn_NN12}
\end{equation}
which measures the discrepancy between the observed data and the model’s predictions. 

Backward pass requires computing $\nabla \mathcal{L}(p)$.
Training a UDE becomes minimizing a cost function $\mathcal{L}(p)$ with respect to $p$.
Similar to standard neural networks, backpropagation applies the chain rule. The main difference is that one must caluculate $d\hat{y}/dp_k$ for $k=1,\dots,s$ , which account for how $\hat{y}(t,p)$ depends on $p$.
By the chain rule, the partial derivatives of the cost with respect to weights and biases take the same rule as in \eqref{eqn_NN8}, except that the last-layer error now includes $\frac{\partial \mathcal{L}}{\partial \hat{y}} \,.\,  \frac{\partial \hat{y}}{\partial u}$ instead of $\frac{\partial \mathcal{L}}{\partial a^{(L)}} $. 
As shown in \eqref{eqn_NN13}, the error at layer $L$ becomes
\begin{equation}
\delta^{(L)} = \frac{\partial \mathcal{L}}{\partial z^{(L)}} = \frac{\partial \mathcal{L}}{\partial \hat{y}} \,.\,  \frac{\partial \hat{y}}{\partial u} \, . \, \sigma'(z^{(L)})
    \label{eqn_NN13}
\end{equation}
and earlier layers follow \eqref{eqn_NN7} as 
\begin{equation}
 \delta^{(l)} = D^{(l)}\left( W^{(l+1)} \right)^T D^{(l+1)}\left( W^{(l+2)} \right)^T \cdots D^{(L-1)}\left( W^{(L)} \right)^T D^{(L)}\left( \frac{\partial \mathcal{L}}{\partial \hat{y}} .  \frac{\partial \hat{y}}{\partial u} \right)^T.
    \label{eqn_NN14}
\end{equation}
Repeated forward evaluation of the UDE, cost calculation, and backward propagation continue until the network meets some convergence threshold or exhausts a computation budget. This hybrid approach allows the learned component of $u$ to fit the training data while respecting known dynamical structure.

\section{Optimal control problem solved by deep neural network}\label{proofOCasDNN}
In this section, we formulate an optimal control problem using the framework of UDEs. Within this framework, the traditional forward-backward sweep method for solving optimal control problems is interpreted through the training of deep neural networks. Specifically, the forward pass solves the state equations, while the backward pass solves the adjoint equations. By optimizing the loss function, which is defined by the objective functional, this approach alleviates the numerical instabilities associated with backward integration of the adjoint equations. Moreover, it combines the descriptive power of differential equations for modeling control dynamics with the adaptability of neural networks to approximate complex, non-linear functions effectively.

\subsection{Formulating the optimal control problem as a model of UDEs}
In the optimal control problem with $ n $ state variables and $ m $ control variables,
we typically define the objective as
\begin{equation}
    \begin{array}{rl}
         \underset{u}{\min} \, \mathcal{L}(u) &= \underset{u}{\min} \, \displaystyle \int_{t_0}^{t_1} f(t, y(t), u(t))\, dt, \\[2.0ex]
        \text{subject to} &  y'(t) = g(t, y(t), u(t)), \quad y(t_0) = y_0.
    \end{array}
        \label{eqn_NN13}
\end{equation}
where $ y(t) = (y_1(t), \dots, y_n(t))^T $, $ u(t) = (u_1(t), \dots, u_m(t))^T $, $ y_0 = (y_{10}, \dots, y_{n0})^T $, and $ g(t, y, u) = (g^1(t, y, u), \dots, g^n(t, y, u))^T $. The functions $f $ and $ g^i $ are assumed to be continuously differentiable in all variables.
For the problem involves  a payoff function $ \phi $, is described as $  \underset{u}{\min} \, \mathcal{L}(u) = \underset{u}{\min} \, \int_{t_0}^{t_1} \tilde{f}(t, y(t), u(t))\, dt + \phi(y(t_1))$.
We employ the Dirac delta function, then rewrite the problem in \eqref{eqn_NN13}, where $f(t, y(t), u(t)) = \tilde{f}(t, y(t), u(t)) + \delta(t - t_1) \phi(y(t))$.

Next, we represent the control function $ u(t, y(t)) $ as a neural network. Given an input with the same number of neurons as the number of controls, $a^{(1)} = (t, \dots, t)^T \; \in \; \mathbb{R}^m$, the output at layer $l$ is defined in the same way as in \eqref{eqn_NN1} for $l = 2, \dots, L$. Then, we define the controls as a neural network as
\begin{equation}
u(t, y(t)) = \text{NN}(t, p),
\label{eqn_NN24}
\end{equation}
where $ p $ includes all the elements in the weight matrices $ W^{(l)} $ and bias vectors $ b^{(l)} $ across the network, with $ p \in \mathbb{R}^s $.
Then the optimal control problem in \eqref{eqn_NN13} is represented in the UDE formulation with the cost function defined by the objective functional as 
\begin{equation}
\begin{array}{l}
    \underset{\text{NN}}{\min} \, \mathcal{L}(\text{NN}) = \underset{\text{NN}}{\min} \, \displaystyle \int_{t_0}^{t_1} f(t, y(t), \text{NN}(t, p))\, dt\\[2.0ex]
    y'(t) = g(t, y(t), \text{NN}(t, p)), \quad y(t_0) = y_0
\end{array}
\label{eqn_NN25}
\end{equation}
or
\begin{equation}
\begin{array}{l}
    \underset{p}{\min} \, \mathcal{L}(p) = \underset{p}{\min} \, \displaystyle \int_{t_0}^{t_1} f(t, y(t), p)\, dt\\[2.0ex]
     y'(t) = g(t, y(t), p), \quad y(t_0) = y_0.
\end{array}
\label{eqn_NN26}
\end{equation}
This reformulation integrates the neural network as a control mechanism, where the objective is to optimize the network parameters to minimize the cost function $ \mathcal{L}(p) $ over the specified time interval $[t_0, t_1]$, thereby aligning with the desired control objective. 
Specifically, we solve the state equation in \eqref{eqn_NN26} forward in time to collect information of the state variable $y(t)$. Then  we calculate the gradients  $\nabla \mathcal{L}(p)$ in \eqref{eqn_NN8} for the update of the network parameters $p$ to minimize the control objective   $ \mathcal{L}(p) $.
Theoretically, the training finishes when we find the network parameter set $p^\ast$, such that $ \mathcal{L}(p^\ast) =0$. We thereby obtain the optimal control $u^\ast(t) = NN(t,\,p^\ast)$ and the corresponding state $y^\ast(t)$, which will be proved to satisfy Pontryagin's maximum principle in the next.


\subsection{Gradient decent}

To choose an update for network parameters $\Delta p$, such that the cost can be smaller than the previous step, we need to calculate the gradient $\nabla \mathcal{L}(p) \in \mathbb{R}^s$ in \eqref{eqn_NN3b}.
We will prove that the network parameter set $p^\ast$, such that $\nabla \mathcal{L}(p^\ast)=0$, provides an optimal control $u^\ast = NN(t, p^\ast)$ and the associated state $y^\ast=y(t, p^\ast)$ satisfying PMP.

With $ n $ states in the system, we introduce $ n $ corresponding adjoint variables, one for each state. This is represented by the piecewise differentiable vector-valued function $ \lambda(t) = [\lambda_1(t), \dots, \lambda_n(t)]^T $, where each $ \lambda_i(t) $ corresponds to the state variable $ y_i(t)$.
For each network parameter $p_r$, $r=1,\dots, s$, the gradient of the loss function  $\nabla \mathcal{L}(p)$ in \eqref{eqn_NN26} is
\begin{equation}
\begin{array}{ll}
     \dfrac{\partial \mathcal{L}}{\partial p_r}&  \displaystyle =   \dfrac{\partial}{\partial p_r} \int_{t_0}^{t_1}f(t,y(t, p),p)dt \\[2.0ex]
     &  \displaystyle = \dfrac{\partial}{\partial p_r} \left( \int_{t_0}^{t_1}
     \left[ f(t,y(t, p),p) +\dfrac{d}{dt} (\lambda(t) \circ y(t, p)) \right]dt 
    -\left[ \lambda(t) \circ y(t, p)\right]|_{t_0}^{t_1} \right)\\[3.0ex]
     &  \displaystyle =\int_{t_0}^{t_1}   \left[ \dfrac{\partial}{\partial p_r} 
     f(t,y(t, p),p) + \dfrac{\partial}{\partial p_r} \left(\dfrac{d \lambda(t)}{dt}  \circ y(t, p) 
     +  \lambda(t) \circ \dfrac{d y(t, p)}{dt} \right) \right]dt 
    -  \dfrac{\partial}{\partial p_r}  \left[ \lambda(t) \circ y(t, p)\right]|_{t_0}^{t_1} \\[3.0ex]
    &  \displaystyle = \int_{t_0}^{t_1}   \left[ \nabla f_y \circ \dfrac{\partial y(t, p)}{\partial p_r} + \dfrac{\partial f}{\partial p_r} + \dfrac{d \lambda(t)}{dt}  \circ \dfrac{\partial y(t, p)}{\partial p_r} + \lambda(t)  \circ \dfrac{d g}{d p_r} \right]dt \\[5.0ex]
    &- \lambda(t_1) \circ \dfrac{\partial y(t_1, p)}{\partial p_r} +\lambda(t_0) \circ \dfrac{\partial y(t_0, p)}{\partial p_r},
\end{array}
    \label{eqn_NN27}
\end{equation}
where $\nabla f_y = (f_{y_1}, \dots, f_{y_n})^T$, $ y(t, p)$ is solved from the state equation in \eqref{eqn_NN26}, and  $ \circ $ denotes the dot product.
By representing the controls by neural networks, we have $ g(t, y(t, p), p) = (g^1(t, y(t, p), p), \dots, g^n(t, y(t, p), p))^T$. The last term in the integrand in \eqref{eqn_NN27} becomes
\begin{equation}
   \begin{array}{rl}
    \lambda(t) \circ \dfrac{d g }{d p_r} =&  \lambda(t) \circ \begin{pmatrix} 
         \nabla g^1_y \circ  \dfrac{\partial y(t)}{\partial p_r} + \nabla g^1_{p_r} \\[2.0ex]
         \vdots \\[2.0ex]
         \nabla g^n_y \circ  \dfrac{\partial y(t)}{\partial p_r} + \nabla g^n_{p_r}
     \end{pmatrix}  \\[11.0ex]
    =&  \lambda(t) \circ \left(  \nabla g^k_y \circ  \dfrac{\partial y(t)}{\partial p_r} \right)^T_{k=1,\dots, n}
     + \lambda(t) \circ  \left( \nabla g^k_{p_r} \right)^T_{k=1,\dots, n} \\[5pt]
      =& \displaystyle \sum\limits_{k=1}^n \lambda_k(t)\,  \nabla g^k_y \circ  \dfrac{\partial y(t)}{\partial p_r}  + \lambda(t) \circ  \left( \nabla g^k_{p_r} \right)^T_{k=1,\dots, n} \\[5pt]
     =& \sum\limits_{k=1}^n \lambda_k(t)\, \left(\sum\limits_{j=1}^n \dfrac{\partial g^k}{\partial y_j}  \dfrac{\partial y_j(t)}{\partial p_r}  \right)  + \lambda(t) \circ  \left( \nabla g^k_{p_r} \right)^T_{k=1,\dots, n} \\[5pt]
    =& \displaystyle  \sum\limits_{j=1}^n \left( \sum\limits_{k=1}^n \lambda_k(t)\,\dfrac{\partial g^k}{\partial y_j}  \right) \dfrac{\partial y_j(t)}{\partial p_r}  + \lambda(t) \circ  \left( \nabla g^k_{p_r} \right)^T_{k=1,\dots, n} \\[5pt]
     =&  \left(\lambda(t) \circ \dfrac{\partial g}{\partial y_j}  \right)_{j=1,\dots,n}^T \circ  \dfrac{\partial y(t)}{\partial p_r}  + \lambda(t) \circ  \left( \nabla g^k_{p_r} \right)^T_{k=1,\dots, n}.
     \end{array}
   \label{eqn_NN28}
\end{equation}
Recalling $y(t_0,p) \equiv y(t_0)$, we have $\frac{\partial y(t_0, p)}{\partial p_r}=0$.
Further, denoting $\left( \nabla g^k_{p_r} \right)^T_{k=1,\dots, n} = \frac{\partial g}{\partial p_r}$, Equation \eqref{eqn_NN27} becomes
\begin{equation}
\begin{array}{rl}
     \dfrac{\partial \mathcal{L}}{\partial p_r} =& \displaystyle \int_{t_0}^{t_1}   \left[ \nabla f_y  + \dfrac{d \lambda(t)}{dt} + \left(\lambda(t) \circ \dfrac{\partial g}{\partial y_j}  \right)_{j=1,\dots,n}^T  \right] \circ \dfrac{\partial y(t)}{\partial p_r} dt  \\[11pt]
   + &\displaystyle   \int_{t_0}^{t_1} \left[ \dfrac{\partial f}{\partial p_r}  + \lambda(t) \circ  \left( \nabla g^k_{p_r} \right)^T_{k=1,\dots, n} \right] dt - \lambda(t_1) \circ \dfrac{\partial y(t_1, p)}{\partial p_r}.
\end{array}
    \label{eqn_NN20}
\end{equation}

\subsection{The necessary conditions}
Now, if the neural network parameter set $p^\ast$ satisfies $\mathcal{L}(p^\ast) = 0$, and gives the optimal control $u^\ast(t) = NN(t, p^\ast)$ and the corresponding optimal state $y^\ast(t)=y(t, p^\ast)$. 
Then, we choose the adjoint function $\lambda(t)$ to satisfy
\begin{equation}
 \dfrac{d\lambda_i(t)}{dt} = - \dfrac{\partial }{\partial y_i} \left[ f(t, y^\ast(t), p^\ast) + \lambda(t) \circ g(t, y^\ast(t), p^\ast) \right],  \; \text{for}\; i = 1, \dots, n  \quad\text{(adjoint equation)},
    \label{eqn_NN21}
\end{equation}
and the boundary condition
\begin{equation}
 \lambda(t) = 0 \quad \text{(transversality condition)}.
    \label{eqn_NN22}
\end{equation}
Now \eqref{eqn_NN20} reduces to 
\begin{equation}
    0= \int_{t_0}^{t_1} \left[ \dfrac{\partial f}{\partial p_r}  + \lambda(t) \circ  \left( \nabla g^k_{p_r} \right)^T_{k=1,\dots, n} \right] \Bigg\rvert_{y^\ast, p^\ast} dt.
    \label{eqn_NN23}
\end{equation}
By denoting $\left( \nabla g^k_{p_r} \right)^T_{k=1,\dots, n}$ as $\frac{\partial g}{\partial p_r}$. Equation $\eqref{eqn_NN23}$ implies the optimality condition
\begin{equation}
  \dfrac{\partial }{\partial p_r}f(t, y^\ast(t), p^\ast)  + \lambda(t)  \circ \dfrac{\partial }{\partial p_r}g(t, y^\ast(t), p^\ast) = 0 \quad \text{for all} \quad t_0 \leq t \leq t_1.
    \label{eqn_NN31}
\end{equation}
The above derivation draws the following results.
\begin{thm}
   If the neural network parameter set $p^\ast$ satisfies $\mathcal{L}(p^\ast) = 0$, which gives the optimal control $u^\ast(t) = NN(t, p^\ast)$ and the corresponding optimal state $y^\ast(t)=y(t, p^\ast)$  for problem \eqref{eqn_NN26}, then there exists a piecewise differentiable adjoint variable $\lambda(t)$ such that
   \begin{equation*}
       H(t,y^\ast(t), u, \lambda(t)) \leq  H(t,y^\ast(t), u^\ast(t), \lambda(t))
   \end{equation*}
   for all controls $u(t) = NN(t, p)$ at each time $t$, where the Hamiltonian $H$ is defined as 
   \begin{equation*}
      H(t, y(t), u(t), \lambda(t)) = f(t, y(t), u(t)) + \lambda(t) \circ g(t, y(t), u(t)),
   \end{equation*}
   and $\lambda(t)$ is governed by the adjoint equation \eqref{eqn_NN21} and satisfies the transversality condition \eqref{eqn_NN22}.
\end{thm}

\section{An application on cancer immunotherapy}\label{Application}
In this section, we develop a UDE model to solve an optimal control problem aimed at determining optimal combination therapy for cancer treatment. The original model exhibits rich dynamical behaviors in immunotherapy, which arise from complex bifurcations [\cite{zhang2024bifurcations}]. Such bifurcations can introduce numerical stiffness, even when solving the state equation forward in time, because the cancer-immune interactions involve highly nonlinear terms. These same nonlinearities also affect the adjoint equation, potentially amplifying instability when that equation is solved backward in time.
To mitigate this instability, we will replace the traditional backward sweep with a backward-propagation approach, using gradient descent to optimize the objective functional directly. This strategy is designed to handle the singularities induced by bifurcations more effectively, ensuring more stable computation.

Cancer immunotherapy activates the immune system, particularly T cells, to identify and destroy cancer cells. The process begins with antigen-presenting cells (APCs) like dendritic cells, which process cancer cell antigens into peptides. These peptides are presented on APC surfaces via major histocompatibility complex (MHC) molecules, enabling T cell recognition.
T cell activation requires two signals: antigen recognition and co-stimulation. The first signal is through T cell receptors interacting with peptides on APCs, while the second signal, co-stimulation, is provided by B7 molecules on APCs interacting with CD28 receptors on T cells. Without both signals, T cells can become anergic and fail to respond effectively.
Cancer-specific T cells then target tumor sites, recognizing and destroying cancer cells. Cancer immunotherapy strategies enhance these responses by promoting co-stimulation and inhibiting co-suppressive signals, like immune checkpoints, to prevent tumors from evading the immune response. Agents used include CAR T cell therapies \cite{harrison2021enhancing}, monoclonal antibodies \cite{effer2023therapeutic}, and immune checkpoint inhibitors like anti-PD-1 and anti-CTLA-4 \cite{curran2010pd}.

Combination therapies aim to simultaneously enhance co-stimulation and inhibit co-suppression, optimizing immune system effectiveness against cancer. However, these therapies can lead to significant side effects, such as cytokine release syndrome from overactive immune responses \cite{maude2014managing} or immune-related adverse events from checkpoint inhibitors \cite{marshall2018introduction}.
Patient responses to these therapies can vary widely, from significant tumor shrinkage to resistance \cite{zheng2021mathematical,zhang2024bifurcations}. Treatment optimization involves balancing immune activation with side effect management, necessitating personalized strategies to maximize benefits while minimizing harm \cite{pappalardo2019silico,alfonso2020translational,jenner2021silico,zahid2021dynamics,bradshaw2019applications,xiong2021translational}.

Mathematical modeling provides a valuable tool for understanding and simulating the complex interactions between immune cells and cancer in the study of immunotherapy. These models allow researchers to explore treatment strategies, optimize dosing, and predict outcomes more efficiently than experimental methods alone. By providing insights into tumor-immune dynamics and guiding the design of personalized therapies, mathematical models help improve the effectiveness and precision of cancer immunotherapy treatments.

\subsection{The model for combination of cancer immunotherapy}
An established immune-cancer model \cite{zheng2021mathematical} is written as follows:
\begin{equation}
    \begin{array}{ll}
        \dfrac{dC}{dt} &= F(C)\,C-\kappa CE,  \\[2.0ex]
        \dfrac{dA}{dt} &= r_A C-\delta_A A, \\[2.0ex]
        \dfrac{dI}{dt} &= r_I CE-\delta_I I, \\[2.0ex]
        \dfrac{dE}{dt} &= -r_E\left(E-E^{\ast} \right) + u_1(t)\beta AIES -u_2(t)\gamma ES,\\[2.0ex]
        \dfrac{dS}{dt} &= -r_S\left(S-S^{\ast} \right) - u_1(t)\beta AIES +u_2(t)\gamma ES.
    \end{array}
    \label{eqn1}
\end{equation}
The mathematical model describes the dynamics of cancer cells $C(t)$, tumor-specific antigen $A(t)$, inflammatory signals $I(t)$, and effector $E(t)$ and non-effector T cells $S(t)$ within the tumor microenvironment. The concentration of cancer cells, $ C $, follows logistic growth up to a maximum rate $ r_{\text{max}} $, which is described as
\begin{equation*}
    F(C)=\min \left\{ r_{max},\, r_C \left(1-\dfrac{C}{C^{\ast}} \right) \right\}.
\end{equation*}
Cancer cells are reduced by the killing effect of effector T cells with a mass-action term proportional to $ \kappa C E $. The antigen concentration, $ A $, produced by antigen-presenting cells, enters the system at a rate proportional to the cancer population, with a production rate $ r_A$, and degrades at a rate $ \delta_A $. Inflammatory signals, $ I $, are produced at the tumor site when cancer cells are killed by effector T cells, modeled in the term $ r_I C E $, and degrade at a rate $ \delta_I $. The effector T cell population, $ E $, tends toward a steady-state level $ E^* $ at a rate $ r_E $ and can increase as non-effector T cells, $ S $, transition to the effector state in the presence of antigen, inflammation, and effector T cells, described by the term $ \beta A I E S $. Conversely, effector T cells can transition back to the non-effector state at a rate proportional to $ \gamma S E $, describing the immune suppression that occurs within the tumor microenvironment. Similarly, non-effector T cells, $ S $, tend toward their steady-state level $ S^* $ at a rate $ r_S $, transitioning to the effector state in the same term $ \beta A I E S $, and transitioning back to the non-effector state at rate $ \gamma$. This model provides a framework to explore the balance between immune activation and suppression in cancer immunotherapy.

Cancer immunotherapy aims to enhance the body's immune response against tumors by modulating two key mechanisms: co-stimulation and co-suppression. Co-stimulation, driven by the parameter $ \beta $ in mathematical models, represents the activation of immune cells, particularly T cells, by signals such as antigens and inflammatory cues in the tumor microenvironment. Increasing $ \beta $ enhances the transition of non-effector (tolerogenic) T cells into effector (immunogenic) T cells, boosting the immune system's ability to recognize and attack cancer cells. On the other hand, co-suppression, controlled by the parameter $ \gamma $, reflects the inhibitory signals that prevent effector T cells from effectively attacking cancer cells. High levels of co-suppression can lead to immune tolerance, where cancer cells evade destruction. Immunotherapies that block co-suppression, such as immune checkpoint inhibitors, effectively reduce the value of $ \gamma $, decreasing the rate at which effector T cells revert to a tolerogenic state. By fine-tuning these parameters, cancer immunotherapy can promote a stronger and more sustained immune response, leading to improved tumor clearance.

\subsection{Optimal treatment strategies for combination of immunotherapy}
Our goal is to develop an optimal strategy for immunotherapeutic treatment. Specifically, we aim to address the following objectives:

\begin{itemize}
    \item To reduce the tumor mass by the end of the treatment period while keeping it under control throughout the treatment duration.
    \item To maximize the population of effector T cells and minimize the population of non-effector T cells both during and at the end of the treatment, with the aim of maintaining remission after the cessation of treatment.
    \item To minimize the treatment burden on patients, taking into account the toxicity of immunotherapy in the design of the method.
\end{itemize}

To achieve these objectives, we formulate an optimal control problem that incorporates a cost function, which includes terms accounting for the final tumor mass, the effector and non-effector T cell populations, and treatment toxicity, both during and at the end of the treatment period. This leads to the following optimal control problem:

\begin{equation}
\min_{u(\cdot) \in \mathcal{U}} \int_0^{t_f} L(x(t,u),u(t))\, dt + \phi(x(T,u)),
    \label{eqn2}
\end{equation}
where $\mathcal{U}$ represents the set of admissible controls with $u(t) = (u_1(t), u_2(t))$, 
$x(t,u) = (C(t,u(t)), A(t,u(t)), \\I(t,u(t)), E(t,u(t)), S(t,u(t)))$ is the solution to Eqn. \eqref{eqn1} corresponding to the control $u(t)$, $L$ represents the cost during the treatment period, $\phi$ is the final cost at the end of the treatment, and $t_f$ is the duration of the treatment.

We focus on patients who do not respond to the current treatment and aim to redesign the treatment strategy using optimal control. To this end, we adjust the effects of the co-stimulation and co-suppression drugs by scaling them to $ u_1 $ and $ u_2 $ times their current effects, respectively. 
The ranges of $ u_1(t) $ and $ u_2(t) $ are:
\begin{equation}
m_1 < u_1(t) < M_1, \quad m_2 < u_2(t) < M_2.
\label{eqn46}
\end{equation}
The doses of the co-stimulation and co-suppression drugs are proportional to
\begin{equation}
v_1(t) = u_1(t) - m_1, \quad v_2(t) = M_2 - u_2(t).
\label{eqn47}
\end{equation}
The  toxicity levels are proportional to $v_1(t)$ and $v_2(t)$.
For co-stimulation, $ u_1 = m_1 $ represents the effect without treatment, while $ u_1 = M_1 $ indicates that the co-stimulation effect has been amplified to $M_1$ times the effect of the current dose. We assume the toxicity of the co-stimulation drug is proportional to its dose then to its effect, i.e., $ v_1(t) = u_1(t)  - m_1$.
For co-suppression, $ u_2 = M_2 $ represents the effect of no treatment, $ u_2 = 0 $ indicates that the increased drug dose halts the effector cell dysfunction process, and $ u_2 = m_2 $ signifies that a further increase in the drug dose reverses the function of non-effector cells, converting them into effector cells.  Thus, the toxicity of the co-suppression drug is proportional to its effects, i.e., $ v_2(t) = M_2 - u_2(t) $.

For the therapeutic effects at the end of the treatment, we define the finite final time $t_f$ and the outcome function as:
\begin{equation}
\phi(x(t_f, u)) = d_1 C(t_f) - d_2 E(t_f) + d_3 S(t_f),
\label{eqn48}
\end{equation}
where the goal is to minimize the cancer cell population and non-effector T cell population while maximizing the effector T cell population. The parameters $ d_1 $, $ d_2 $, and $ d_3 $ are weights that reflect the relative importance of these objectives at the end of the treatment period.

For therapeutic effects during the treatment, the objective remains to minimize the tumor mass and non-effector T cell population while maximizing effector T cells. We define the cost function for this purpose as:
\begin{equation}
L_1(x(t, u)) = a C(t) - b E(t) + c S(t),
\label{eqn49}
\end{equation}
where $ a $, $ b $, and $ c $ are weight parameters that balance the different aspects of the treatment objectives.

To account for toxicity penalties, one option is to use a quadratic penalty for drug dosage, which results in the following cost function:
\begin{equation}
L_2(u(t)) = c_1 v_1(t)^2 + c_2 v_2(t)^2,
\label{eqn50}
\end{equation}
where $ c_1 $ and $ c_2 $ are weight parameters reflecting the toxicity levels of  the co-stimulation and co-suppression drugs, respectively. Quadratic terms are often preferred due to their convexity, which simplifies mathematical and numerical solutions to the optimal control problem.
Considering the combination of therapeutic efficacy and toxicity penalties during the treatment and at the end of the treatment, we define the objective functional as:
\begin{equation}
J_1(u_1, u_2) = 
    \int_0^{t_f} \left(L_1 + L_2 \right) dt + \phi(x(t_f, u)),
\label{eqn51}
\end{equation}
which is to be minimized over all $ (u_1, u_2) \in U $.
Here, the control set $ U $ is defined as:
\begin{equation}
U = \left\{(u_1(t), u_2(t)) \mid m_1 \leq u_1(t) \leq M_1, \, m_2 \leq u_2(t) \leq M_2 \right\},
\label{eqn52}
\end{equation}
A summary of our baseline parameter estimates is presented in Table \ref{tab:baseline}. All simulations in our study are based on these values. Additionally, Table \ref{tab:initial} provides a summary of the initial values used in the simulations.

We choose parameter set which leads to no response for model \eqref{eqn1} with no control.
The control is formed by a neural network after a sliding transformation to meet the requirement in \eqref{eqn46}.
The neural network consists of an input layer, two hidden layers with $10$ neurons and output layer. 
The Gaussian error linear unit (GELU) have been used as the activation function in the first two hidden layers due to its high-performance advantage \cite{hendrycks2016gaussian}.
For the last hidden layer, we choose a tanh activation function to restrict the range of the output of the neural network in $(-1,\,1)$.
Finally, a sliding transformation is applied on the output of the neural network as follows:
\begin{equation}
\begin{array}{l}
    NN_p(\cdot, t) = \text{tanh}\left(
    W^{(3)} \, 
    \text{GELU} \left( 
    W^{(2)} \, \text{GELU} \left( W^{(1 )}
    \begin{bmatrix}
        t\\t
    \end{bmatrix} + b^{(1)} \right)
    + b^{(2)} \right) + b^{(3)} \right)  \\[2.0ex]
   u_1(t) = u_1(NN_p[1,t]) = m_1 + \dfrac{M_1-m_1}{2} \left(NN_p[1,t]  + 1 \right), \\[2.0ex]
u_2(t) = u_2(NN_p[2,t]) = M_1 + \dfrac{M_1-m_1}{2} \left(NN_p[2,t]  - 1 \right). 
\end{array}
\label{eqn53}
\end{equation}
Here, the neural network $ NN_p $ is generated using \texttt{Lux.jl} \cite{pal2023lux,pal2023efficient} with a parameter set $ p $, which consists of elements in the weights and biases. We have $ W^{(1)}\in \mathbb{R}^{10 \times 2}$, $ W^{(2)} \in \mathbb{R}^{10 \times 10}$,   $ W^{(3)} \in \mathbb{R}^{2 \times 10}$, $b^{(1)} \in \mathbb{R}^{10}$, $b^{(2)}\in \mathbb{R}^{10}$, and $b^{(3)}\in \mathbb{R}^2$.
We started by using the discrete form of  the loss function in \eqref{eqn51}  with $l_2$ norm. We optimize the parameters in the set $p$ to minimize the loss function in \eqref{eqn51}  using Adam with a learning rate of $0.01$ for 100 iterations and BFGS with an initial step norm of $0.01$ until the result of the loss function is stabilized.

\subsection{Immunotherapy strategy inducing responses}
\begin{figure}[ht]
    \centering
\includegraphics[width=0.32\textwidth]{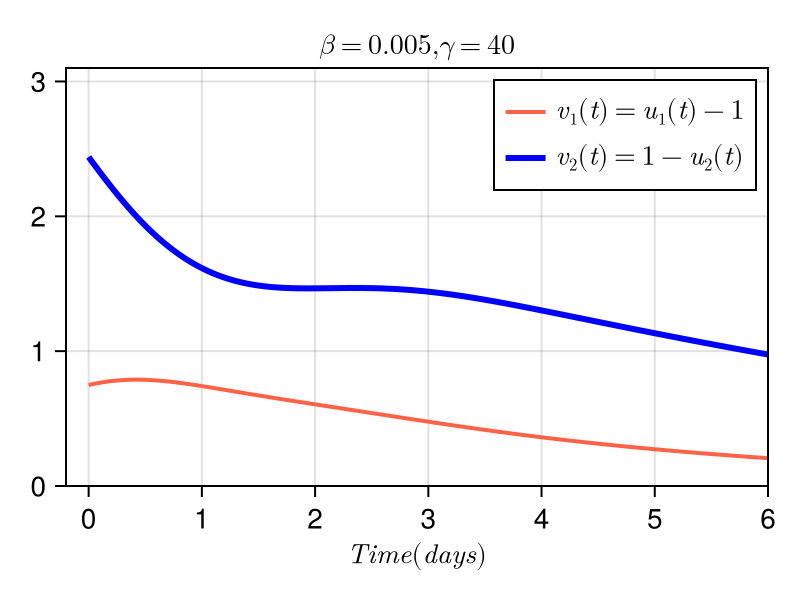}
\includegraphics[width=0.32\textwidth]{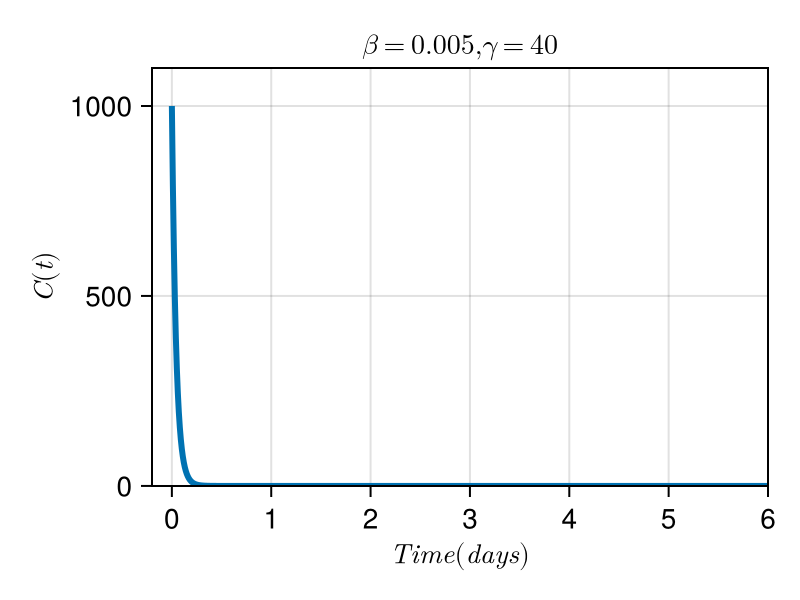}
\includegraphics[width=0.32\textwidth]{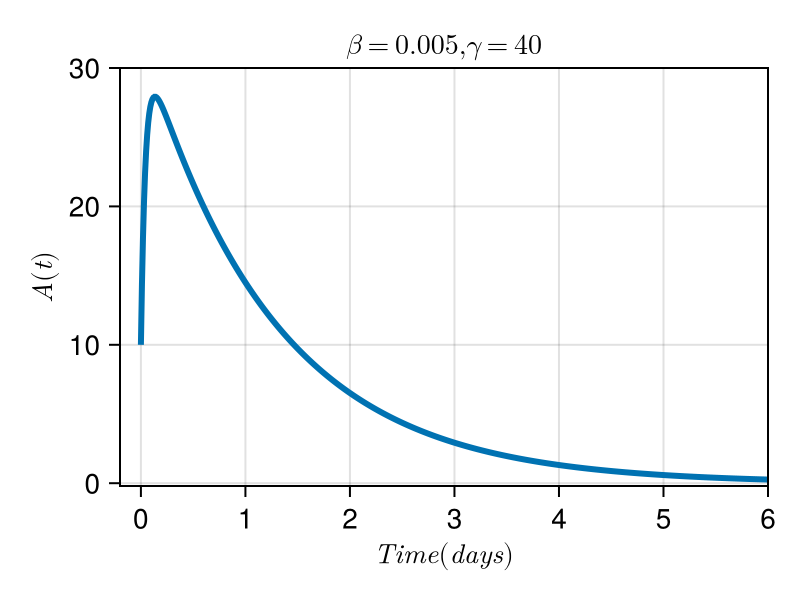}
\includegraphics[width=0.32\textwidth]{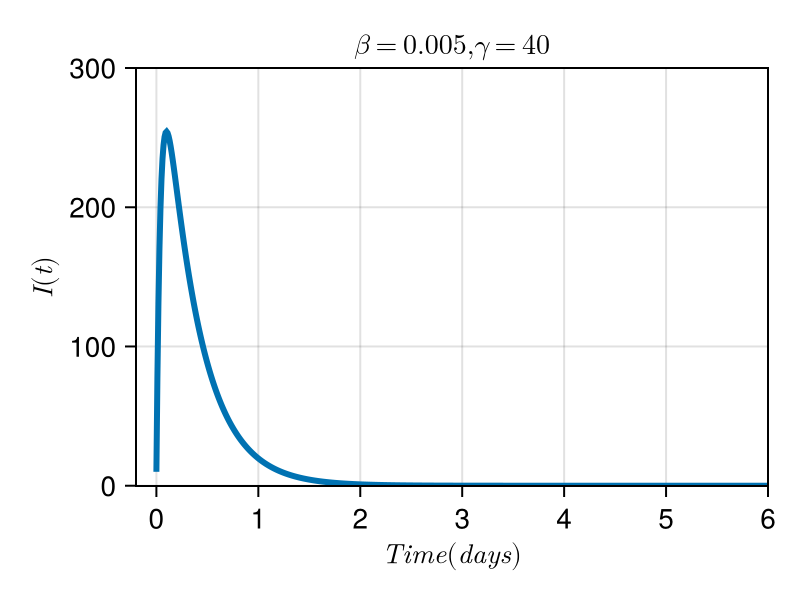}
\includegraphics[width=0.32\textwidth]{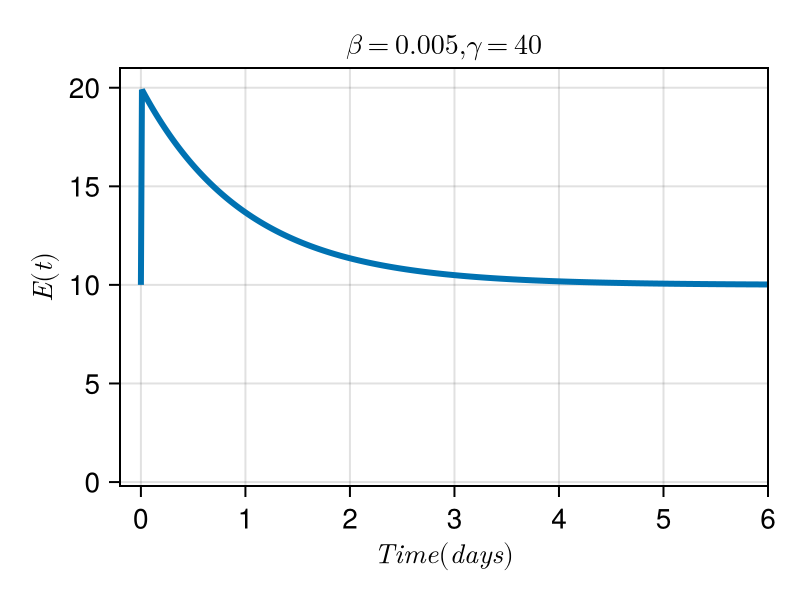}
\includegraphics[width=0.32\textwidth]{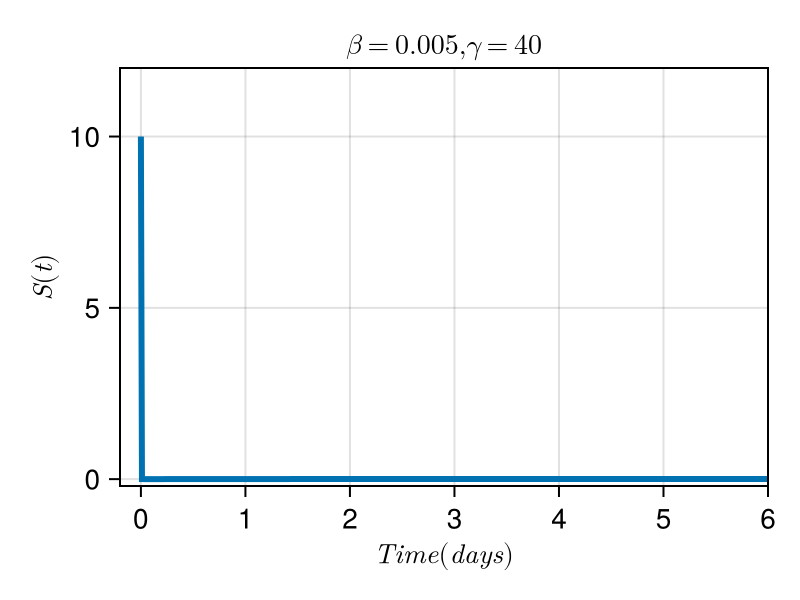}
    \caption{Optimal control strategy under a quadratic cost penalty with moderate-to-high dosing range for co-stimulation and co-suppression therapy. The control successfully reduces cancer cell and non-effector T cell populations while enhancing effector T cells throughout the treatment period. This illustrates an effective immunotherapy response under optimized conditions. Optimal control with quadratic penalty in \eqref{eqn51}, where $a=1$, $b = 10$, $c=100$, $c_{1}=2$, $c_{2}=1$, $d_1=1$, $d_2=10$, $d_3=100$. }
    \label{fig1}
\end{figure}

In our numerical experiments, we explore the effects of immunotherapy, which exhibits a wide range of treatment effects. 
In general, toxicities with  cosupression drug, such as anti-PD-1/PD-L1 mAbs appear to be less common and less severe when compared with costimulation drugs, such as anti-CTLA-4 mAbs \cite{naidoo2015toxicities}.
We set the dosage range for the co-stimulation drug between $m_1=1$ and $M_1=3$, and for the co-suppression drug between $m_2=-3$ and $M_2=1$. We model treatments, that are administered every other day \cite{wu2018immunogenic}, as continuous control. As depicted in Figures \ref{fig1}, our findings indicate that administering the maximum dose of the drugs can rapidly reduce the populations of cancer cells and non-effector cells while enhancing the population of effector cells.

\subsection{Immunotherapy strategy inducing no response}
\begin{figure}[ht]
    \centering
\includegraphics[width=0.32\textwidth]{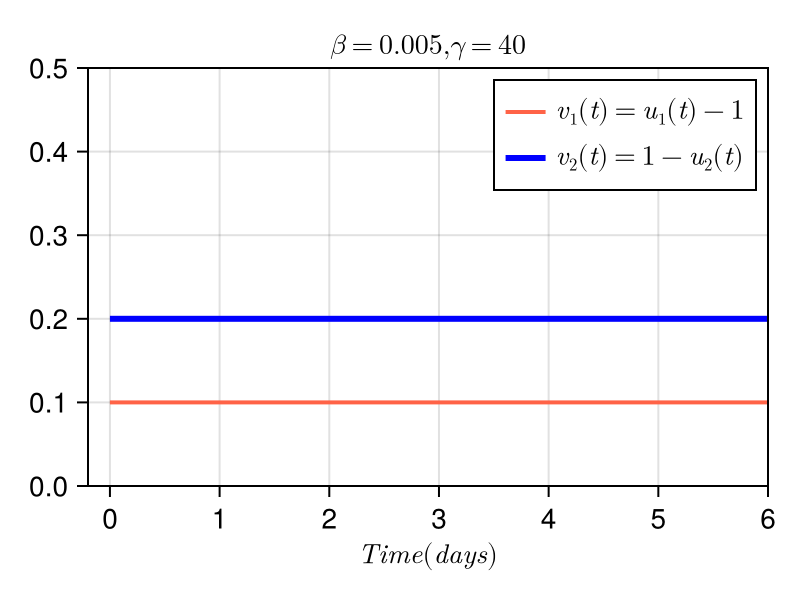}
\includegraphics[width=0.32\textwidth]{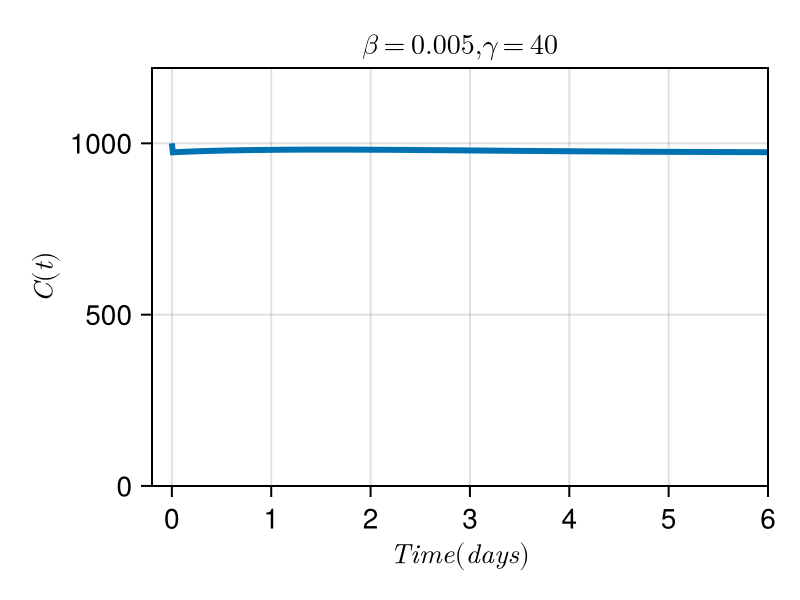}
\includegraphics[width=0.32\textwidth]{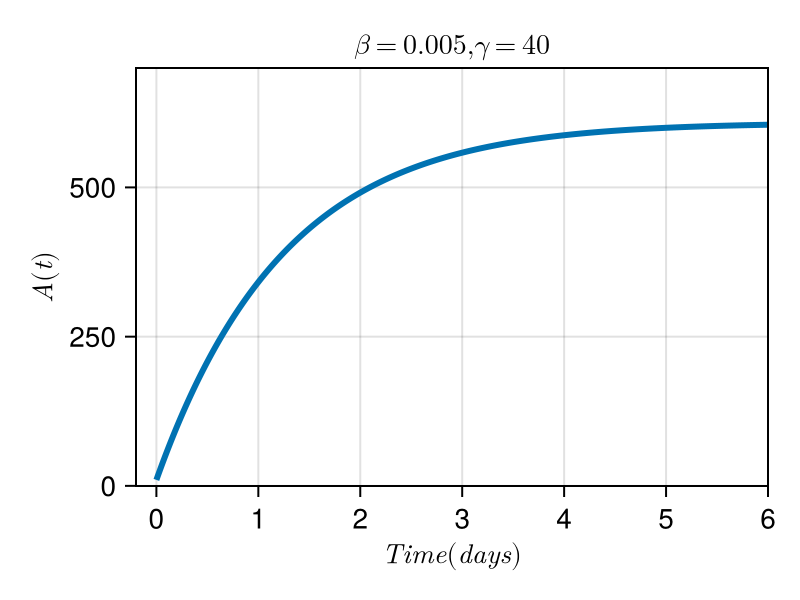}
\includegraphics[width=0.32\textwidth]{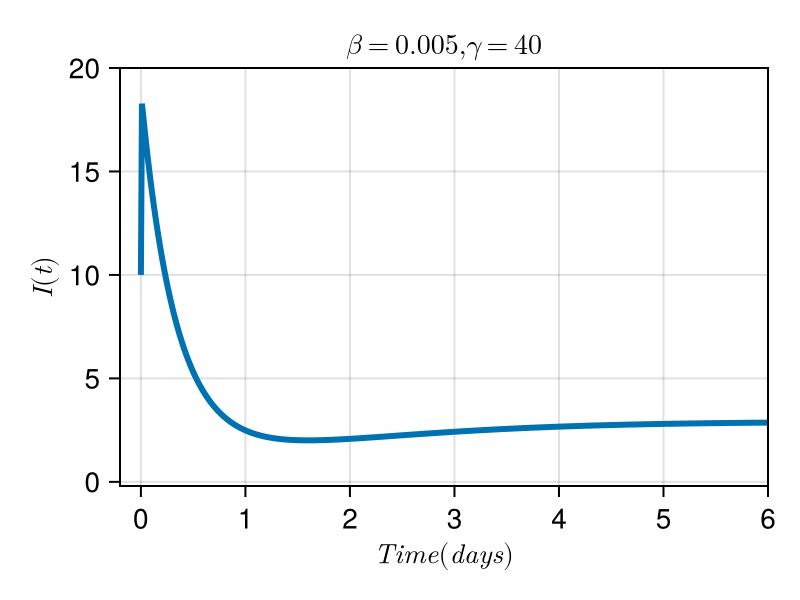}
\includegraphics[width=0.32\textwidth]{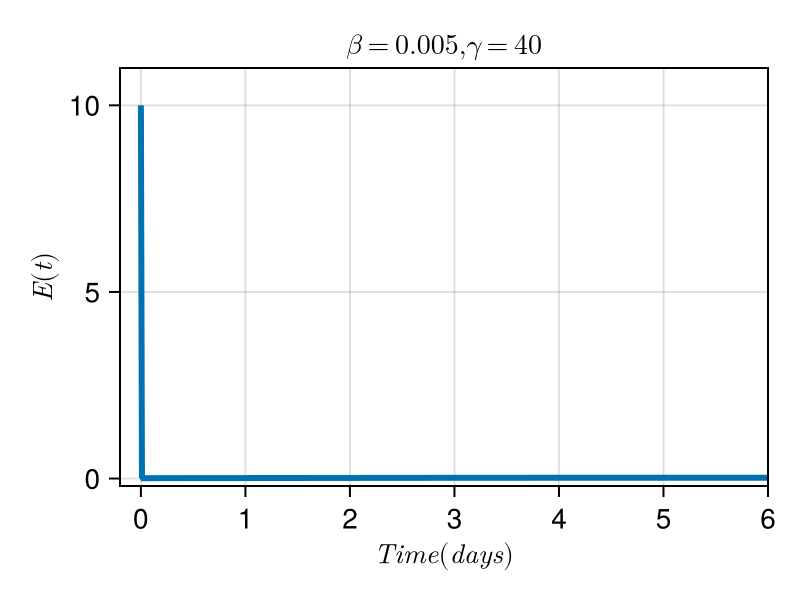}
\includegraphics[width=0.32\textwidth]{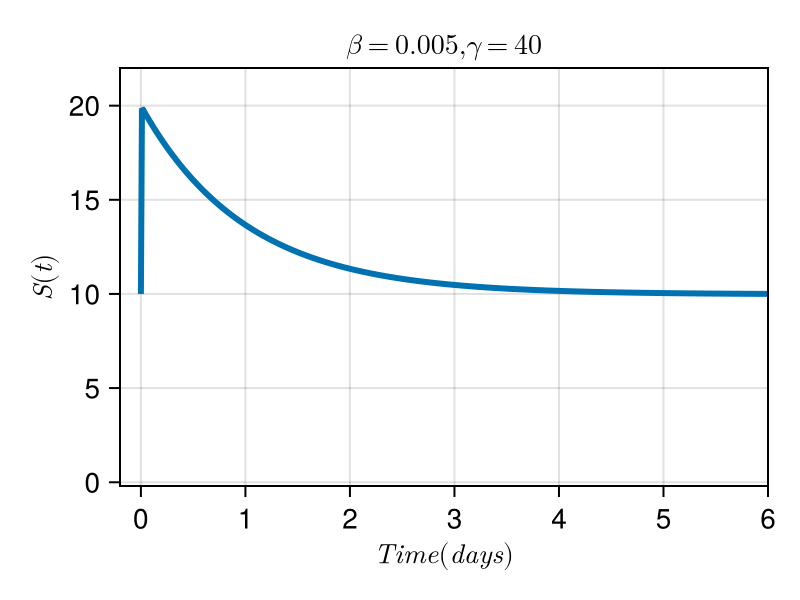}
    \caption{Suboptimal control using low-dose immunotherapy. The limited dosing range fails to control cancer growth, resulting in persistent tumor burden, elevated antigen and non-effector T cell levels, and poor effector T cell activation. This scenario reflects a non-responder case despite treatment. Optimal control with quadratic penalty in \eqref{eqn51}, where $a=1$, $b = 10$, $c=100$, $c_{1}=2$, $c_{2}=1$, $d_1=1$, $d_2=10$, $d_3=100$. }
    \label{fig2}
\end{figure}

Given the relative ease of administration and the lower toxicity, low-dose drugs are often used for maintenance or consolidation therapy \cite{wu2018immunogenic}. We thus consider narrower ranges of treatment effects: $m_1 = 1$ and $M_1 = 1.1$ for the co-stimulation drug, and $m_2 = 0.8$ and $M_2 = 1$ for the co-suppression drug. The numerical experiment depicted in Figure \ref{fig2} illustrates that the low-dose therapy fails to control the cancer cell population, leading to elevated levels of antigens and inflammation signals. In this scenario, the effector cell level remains low while the non-effector cell level is high.

\subsection{Combination of chemotherapy with immunotherapy}

In various clinical studies, low doses of chemotherapeutic agents have been used to precondition patients for immunotherapy. Notably, a study involving 28 metastatic breast cancer patients utilized a combination of low-dose cyclophosphamide with a HER2-positive, allogeneic, GM-CSF-secreting tumor vaccine, demonstrating enhanced responses. Conversely, higher doses of cyclophosphamide were found to suppress these responses \cite{emens2009timed}.  
Building on these findings, we designed a numerical experiment to investigate the combined effects of chemotherapy with co-stimulation and co-suppression immunotherapy. The effects of co-stimulation and co-suppression in our model are represented by $u_1(t)$ and $u_2(t)$, respectively, as defined in \eqref{eqn1}. Additionally, the influence of chemotherapy is incorporated as $u_3(t)$ in model \eqref{eqn54}.
Some traditional chemotherapy drugs, like doxorubicin and cyclophosphamide, have the ability to boost the immunogenic response of cancer cells through a process known as immunogenic cell death (ICD) \cite{wu2018immunogenic}. This process involves several key steps: calreticulin moves to the surface of tumor cells, signaling dendritic cells for tumor antigen uptake; ATP is released, promoting the activity of macrophages \cite{elliott2009nucleotides} and natural killer cells and encouraging IFN$\gamma$ production \cite{beavis2012cd73}; and protein HMGB1 (high mobility group box 1) is expelled post-apoptosis, activating various toll-like receptors and enhancing immune responses against cancer \cite{yanai2009hmgb}. In addition, these drugs can enhance type I interferon pathways, helping to increase the general activation of antitumor immunity in the body \cite{li2013hmgb1}.
The effect of chemotherapy on the activity of macrophages and boosting antigen presenting cells is denoted as $u_A(t)$, on the activation of inflammation signal is denoted as $u_I(t)$. 
The model for combination pf chemotherapy with immunotherapy is written as
\begin{equation}
    \begin{array}{ll}
        \dfrac{dC}{dt} &= u_3(t)\,F(C)\,C-\kappa CE,  \\[2.0ex]
        \dfrac{dA}{dt} &= u_A(t)\,r_A C-\delta_A A, \\[2.0ex]
        \dfrac{dI}{dt} &= u_I(t)\,r_I CE-\delta_I I, \\[2.0ex]
        \dfrac{dE}{dt} &= -r_E\left(E-E^{\ast} \right) + u_1(t) \beta AIES - u_2(t) \gamma ES,\\[2.0ex]
        \dfrac{dS}{dt} &= -r_S\left(S-S^{\ast} \right) - u_1(t)  \beta AIES + u_2(t)  \gamma ES,
    \end{array}
    \label{eqn54}
\end{equation}
where
\begin{equation}
\centering
\begin{array}{cl}
    u_1(t) &= m_1 + \dfrac{M_1-m_1}{2} \left(NN[1,t]  + 1 \right), \\[2.0ex]
    u_2(t) &= M_2 + \dfrac{M_2-m_2}{2} \left(NN[2,t]  - 1 \right), \\[2.0ex] 
    u_3(t) &= m_3 + \dfrac{M_3 - m_3}{2} \, (NN_p[3,t] + 1), \\[2.0ex]
    u_A(t) &=   1 + e_1 \, \text{tanh}(-\text{ln}(u_3(NN_p[3,t]))), \\[2.0ex]
    u_I(t) &=  1 + e_2 \, \text{tanh}(-\text{ln}(u_3(NN_p[3,t]))).
    \end{array}
    \label{eqn55}
\end{equation}
Here, the control set $ U $ is defined as:
\begin{equation}
U_2 = \left\{(u_1(t), u_2(t), u_3(t)) \mid m_1 \leq u_1(t) \leq M_1, \, m_2 \leq u_2(t) \leq M_2 , \, m_3 \leq u_3(t) \leq M_3 \right\}.
\label{eqn52b}
\end{equation}
The doses of the co-stimulation and co-suppression drugs and chemotherapy are proportional to $v_1(t)$ and $v_2(t)$ in \eqref{eqn47} and $v_3(t)= M_3-u_3(t)$.
We define the objective functional as
\begin{equation}
J_2(u_1, u_2, u_3) = 
      \int_0^{t_f} \left(L_1 + L_3 \right) dt + \phi(x(t_f, u)),
      \label{eqn56}
\end{equation}
where $\phi$ and $L_1$ are taken in \eqref{eqn48} and \eqref{eqn49} and $L_3(u(t)) = c_1 v_1(t)^2 + c_2 v_2(t)^2 + c_3 v_3(t)^2$. $J_2$ is to be minimized over all $ (u_1, u_2, u_3) \in U_2 $.
We simulate low-dose drug treatment combination by taking $m_1=1$, $M_1=1.1$, $m_2=-0.8$, $M_2=1$, $m_3=0.7$, and $M_3=1.1$.
Our result shows that cancer cell can be minimized, however the effector cell population can not be enhanced.
\begin{figure}[ht]
    \centering
\includegraphics[width=0.32\textwidth]{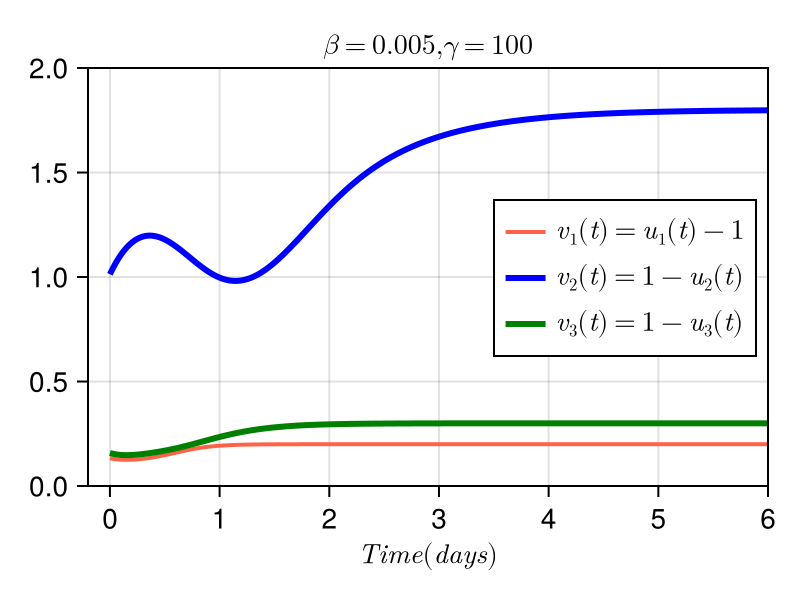}
\includegraphics[width=0.32\textwidth]{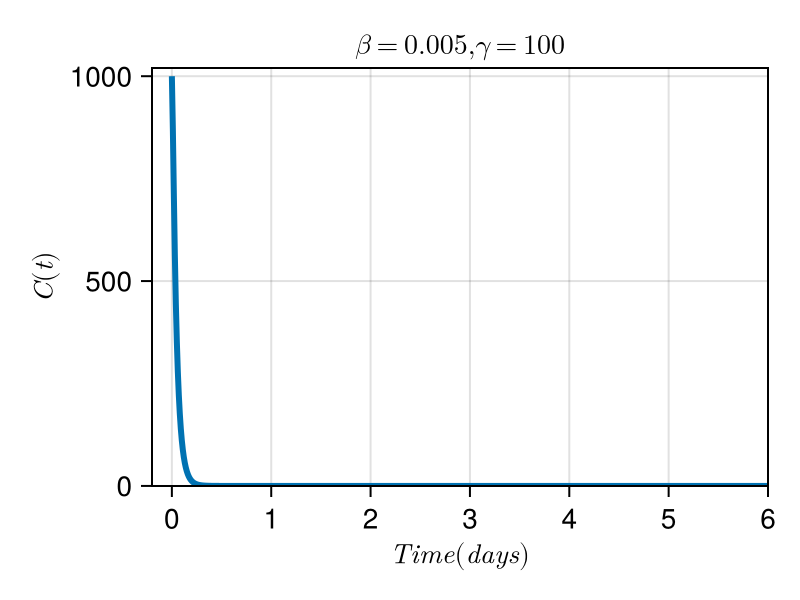}
\includegraphics[width=0.32\textwidth]{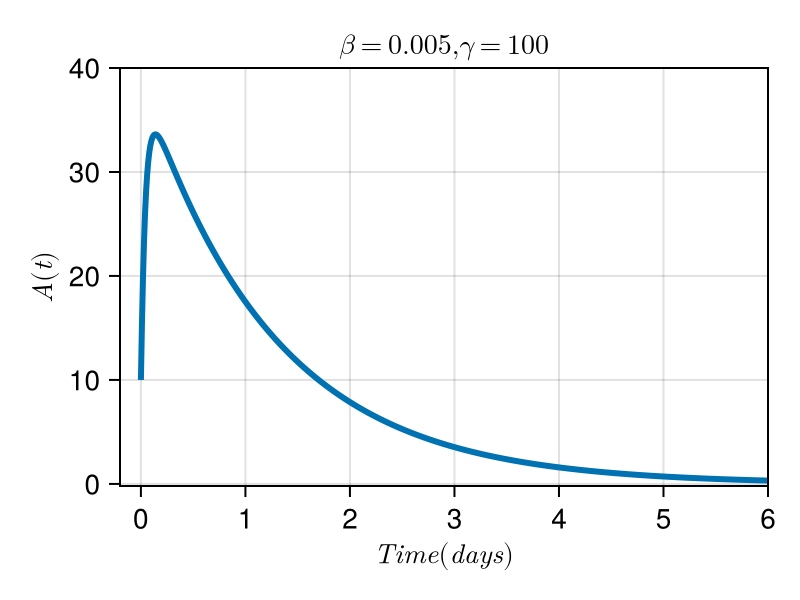}
\includegraphics[width=0.32\textwidth]{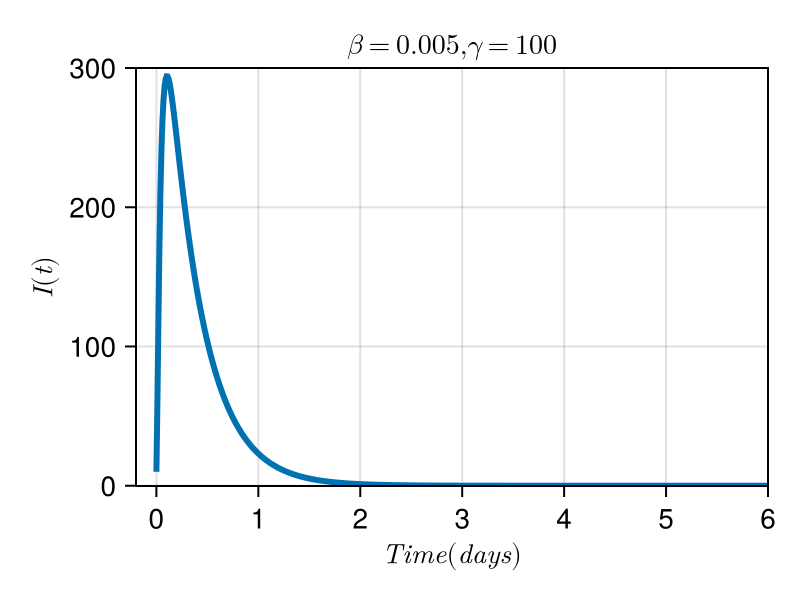}
\includegraphics[width=0.32\textwidth]{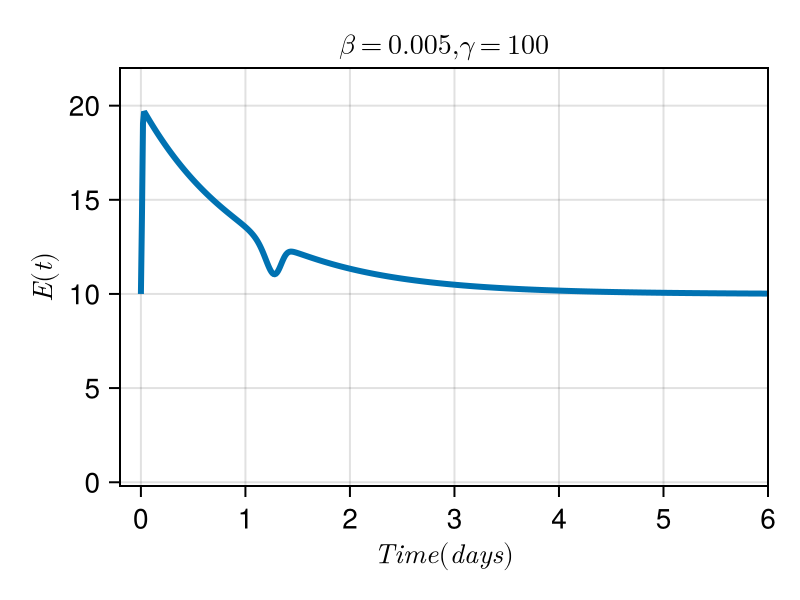}
\includegraphics[width=0.32\textwidth]{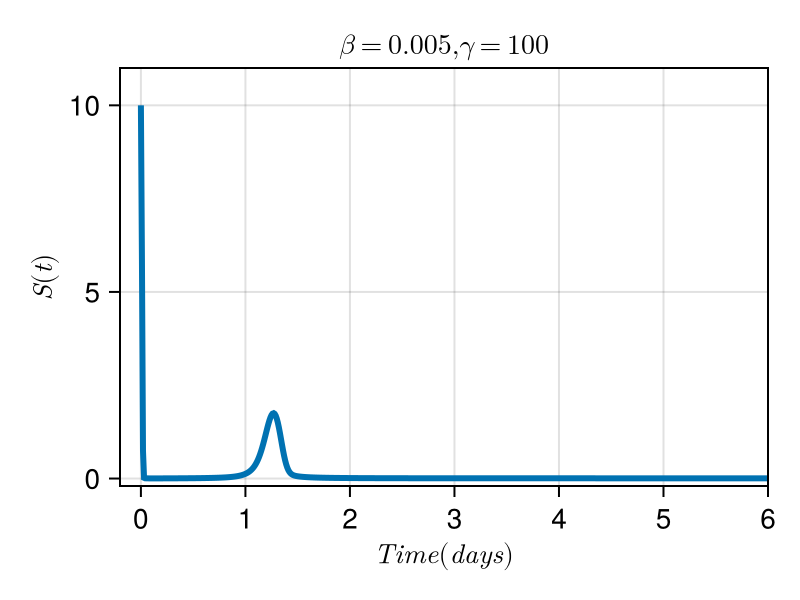}

    \caption{Optimal combination of low-dose chemotherapy with immunotherapy. The control strategy effectively reduces cancer cell burden and substantially enhances the effector T cell population. Here $a=1$, $b = 1$, $c=100$, $c_{1}=2$, $c_{2}=1$, $c_{3}=1$,$e_{1}=2$, $e_{2}=1$, $d_1=1$, $d_2=1$, $d_3=100$. }
    \label{fig3}
\end{figure}

\section{Conclusion}\label{concl}
Universal differential equations (UDEs) integrate a universal function approximator, such as a neural network, within a differential equation framework. In this project, we frame the optimal differential equations (UDEs) integrate a universal function approximator, such as a neural network, within a differential equation framework. This project frames the optimal control problem using UDEs and employs neural networks to approximate control variables. We show that the conditions for an optimal solution, derived from Pontryagin's Maximum Principle, can be met by training the neural network to minimize the objective function defined as the loss function in the UDEs.

Furthermore, by formulating optimal control problems within the UDE framework, we can apply established numerical algorithms to solve the differential equations. This enables us to address state constraints using numerical solvers that handle various types of differential equations, including ordinary, delay, partial, stochastic, and discontinuous equations.
We implement this deep learning approach in the Julia programming language on an immune-cancer dynamics model to design an optimal therapy combination. This novel method efficiently identifies the optimal solution for chemo-immuno combination therapy, affecting multiple factors in cancer-immune dynamics.
This study lays the groundwork for employing deep learning to solve optimal control problems with state constraints in diverse forms.

Beyond cancer therapy, our approach is broadly applicable in fields where optimal control under complex dynamics is essential. It easily adapts to biomedical challenges, optimizing treatments for chronic infections, autoimmune disorders, and metabolic diseases. In public health, it enhances intervention planning, such as vaccination strategies and epidemic containment. Additionally, the method proves valuable in ecology for controlling invasive species and promoting sustainable harvesting, as well as in engineering for ensuring the safety and longevity of lithium-ion battery design. It also finds application in financial systems for risk management and resource optimization, serving as a versatile tool in scientific machine learning for solving high-dimensional, nonlinear, and uncertain dynamic problems.

\begin{table}[ht]
\centering
\caption{Summary of baseline parameter estimates \cite{zheng2021mathematical}.}
\begin{tabular}{|c|l|c|}
\hline
\textbf{Symbol} & \textbf{Description} & \textbf{Estimate} \\ \hline
$r_C$ & Logistic growth rate of cancer & 1.0 day$^{-1}$ \\ \hline
$r_{\text{max}}$ & Maximum growth rate of cancer & 0.09 day$^{-1}$ \\ \hline
$C^*$ & Cancer steady state concentration & $10^3$ cells/nL \\ \hline
$\kappa$ & Killing rate of cancer cells by T cells & 1.2 nL/(cells$\cdot$day) \\ \hline
$r_A$ & Antigen presentation source rate & 0.5 pep/(nL$\cdot$day) \\ \hline
$\delta_A$ & Antigen presentation degradation rate & 0.8 day$^{-1}$ \\ \hline
$r_I$ & Inflammation source rate & 0.4 (ng$\cdot$nL)/(cells$^2\cdot$day) \\ \hline
$\delta_I$ & Inflammation degradation rate & 3.0 day$^{-1}$ \\ \hline
$r_E$ & Effector T cell growth coefficient & 1.0 day$^{-1}$ \\ \hline
$E^*$ & Effector T cell base steady state & 5.0 cells/nL \\ \hline
$r_S$ & Non-effector T cell growth coefficient & 1.0 day$^{-1}$ \\ \hline
$S^*$ & Non-effector T cell base steady state & 5.0 cells/nL \\ \hline
$\beta$ & Effector T cell recruitment coefficient & 0.009 nL$^3$/(pep$\cdot$ng$\cdot$cells$\cdot$day) \\ \hline
$\gamma$ & Non-effector T cell recruitment coefficient & 37.414 nL/(cells$\cdot$day) \\ \hline
\end{tabular}
\label{tab:baseline}
\end{table}

\begin{table}[!ht]
\centering
\caption{Summary of initial values.}
\begin{tabular}{|c|l|}
\hline
\textbf{Variable} & \textbf{Value} \\ \hline
$C$ & 1000  cells/nL \\ \hline
$A$ & 10 peptides/nL  \\ \hline
$I$ & 10 ng/nL  \\ \hline
$E$ & 10 cells/nL \\ \hline
$S$ & 10 cells/nL\\ \hline
\end{tabular}
\label{tab:initial}
\end{table}

    \end{document}